\documentclass[12pt]{iopart}
\usepackage{iopams}

\expandafter\let\csname equation*\endcsname\relax
\expandafter\let\csname endequation*\endcsname\relax
\usepackage{amsmath,amsfonts,amsthm,amsopn,graphicx}
\usepackage{cite}

\newtheorem{definition}{Definition}
\newtheorem{property}{Property}
\newtheorem{theorem}{Theorem}
\newtheorem{lemma}[theorem]{Lemma}

\DeclareMathOperator*{\argmin}{\arg\min}

\usepackage{color}
\usepackage{makecell}

\usepackage[ruled]{algorithm2e}

\begin{document}

\title[Solutions to Least-Squares Problems with Mixed Linear and Antilinear Operators]{Efficient Iterative Solutions to Complex-Valued Nonlinear Least-Squares Problems with Mixed Linear and Antilinear Operators}

\author{Tae~Hyung~Kim and Justin~P.~Haldar}

\address{Ming Hsieh Department of Electrical and Computer Engineering, University of Southern California, Los Angeles, CA, 90089, USA}
\ead{taehyung@usc.edu and jhaldar@usc.edu}
\vspace{10pt}
\begin{indented} 
\item[]April 2020
\end{indented}

\begin{abstract}
We consider a setting in which it is desired to find an optimal complex vector $\mathbf{x}\in\mathbb{C}^N$ that satisfies $\mathcal{A}(\mathbf{x}) \approx \mathbf{b}$ in a least-squares sense, where $\mathbf{b} \in \mathbb{C}^M$ is a data vector (possibly noise-corrupted), and $\mathcal{A}(\cdot): \mathbb{C}^N \rightarrow \mathbb{C}^M$ is a measurement operator.  If $\mathcal{A}(\cdot)$ were linear, this reduces to the classical linear least-squares problem, which has a well-known analytic solution as well as powerful iterative solution algorithms.  However, instead of linear least-squares,  this work considers the more complicated scenario where $\mathcal{A}(\cdot)$ is nonlinear, but can be represented as the summation and/or composition of some operators that are linear and some operators that are antilinear.  Some common nonlinear operations that have this structure include complex conjugation or taking the real-part or imaginary-part of a complex vector.  Previous literature has shown that this kind of mixed linear/antilinear least-squares problem can be mapped into a linear least-squares problem by considering $\mathbf{x}$ as a vector in $\mathbb{R}^{2N}$ instead of $\mathbb{C}^N$.  While this approach is valid, the replacement of the original complex-valued optimization problem with a real-valued optimization problem can be complicated to implement, and can also be associated with increased computational complexity.  In this work, we describe theory and computational methods that enable mixed linear/antilinear least-squares problems to be solved iteratively using standard linear least-squares tools, while retaining all of the complex-valued structure of the original inverse problem.  An illustration is provided to demonstrate that this approach can simplify the implementation and reduce the computational complexity of iterative solution algorithms.

\vspace{2pc}
\end{abstract}

\section{Introduction}

Consider a generic complex-valued finite-dimensional inverse problem scenario in which the forward model is represented as 
\begin{equation}
\mathbf{b} = \mathcal{A}(\mathbf{x}) + \mathbf{n},
\end{equation}
where $\mathbf{b} \in \mathbb{C}^M$ represents the measured data, $\mathcal{A}(\cdot): \mathbb{C}^N \rightarrow \mathbb{C}^M$ is the measurement operator, $\mathbf{n} \in \mathbb{C}^M$ represents noise, and $\mathbf{x} \in \mathbb{C}^N$ represents the unknown signal that we wish to estimate based on knowledge of $\mathbf{b}$ and $\mathcal{A}(\cdot)$.  A common approach to solving this inverse problem is to find a least-squares solution
\begin{equation}
\hat{\mathbf{x}} = \argmin_{\mathbf{x} \in \mathbb{C}^N} \|\mathcal{A}(\mathbf{x}) - \mathbf{b}\|_2^2,\label{eq:ls}
\end{equation}
where $\|\cdot\|_2$ denotes the standard $\ell_2$-norm.  This choice of formulation can be justified in multiple ways, and e.g., corresponds to  the optimal maximum likelihood estimator when the noise vector $\mathbf{n}$ is independent and identically-distributed (i.i.d.) Gaussian noise \cite{kay1993}.  Even for more complicated noise statistics that follow, e.g., the Poisson, Rician, or non-Central Chi distributions, there exist iterative methods that allow the maximum likelihood estimator to be obtained by iteratively solving a sequence of least-squares objective functions \cite{erdogan1999,fessler1998,varadarajan2015}.    In addition, another reason for the popularity of least-squares is that the optimization problem is frequently very easy to solve.  For example, in the case where $\mathcal{A}(\cdot)$ is a linear operator (i.e., $\mathcal{A}(\cdot)$ can be represented in an equivalent matrix form as $\mathcal{A}(\mathbf{x}) = \mathbf{A}\mathbf{x}$ for some matrix $\mathbf{A}\in\mathbb{C}^{M \times N}$) with a trivial nullspace, the solution to Eq.~\eqref{eq:ls} has the analytic closed-form expression \cite{luenberger1969}
\begin{equation}
\hat{\mathbf{x}} = (\mathbf{A}^H\mathbf{A})^{-1}\mathbf{A}^H \mathbf{b}, \label{eq:norm}
\end{equation} 
where $^H$ denotes the conjugate-transpose operation.  In large-scale problems where $N$ is very large, the matrix inversion in Eq.~\eqref{eq:norm} may be computationally intractable, although there exist a variety of simple iterative algorithms that are guaranteed to converge to a globally-optimal solution, including Landweber iteration \cite{landweber1951}, the conjugate gradient (CG) algorithm \cite{hestenes1952}, and LSQR \cite{paige1982}.

Instead of assuming linearity,  we focus in this work on solving least-squares problems in the scenario where $\mathcal{A}(\cdot)$ is nonlinear, but can be represented as the summation and/or composition of some operators that are linear and some operators that are antilinear.   Such nonlinear operators have sometimes been termed as \emph{real-linear} operators in mathematical physics \cite{huhtanen2011}.  Important common examples of operators that possess this kind of nonlinear structure include the complex-conjugation operator  
\begin{equation}
\mathcal{A}(\mathbf{x}) = \overline{\mathbf{x}}, 
\end{equation}
the operator that takes the real part of a complex vector
\begin{equation}
\mathcal{A}(\mathbf{x}) = \mathrm{real}(\mathbf{x}) \triangleq \frac{1}{2}\mathbf{x} +  \frac{1}{2}\overline{\mathbf{x}},\label{eq:re}
\end{equation}
and the operator that takes the imaginary part of a complex vector
\begin{equation}
\mathcal{A}(\mathbf{x}) = \mathrm{imag}(\mathbf{x}) \triangleq  \frac{1}{2i}\mathbf{x} -  \frac{1}{2i}\overline{\mathbf{x}}.\label{eq:im}
\end{equation}

Even though the descriptions we present in this paper are generally applicable to arbitrary real-linear operators, we were initially motivated to consider such operators because of specific applications in magnetic resonance imaging (MRI) reconstruction.  In particular, MRI images are complex-valued, and real-linear operators have previously been used to incorporate prior information about the image phase characteristics into the image reconstruction process, which helps to regularize/stabilize the solution when the inverse problem is ill posed.  For example, there is a line of research within MRI that poses phase-constrained image reconstruction as \cite{bydder2005, willig2005, lew2007, hoge2007, haldar2012,blaimer2016} 
\begin{equation}
\begin{split}
\hat{\mathbf{x}} &= \argmin_{\mathbf{x} \in \mathbb{C}^N} \|\mathbf{A}\mathbf{x} - \mathbf{b}\|_2^2 + \lambda \|\mathrm{imag}(\mathbf{B}\mathbf{x})\|_2^2\\
&= \argmin_{\mathbf{x} \in \mathbb{C}^N} \left\|\begin{bmatrix} \mathbf{A}\mathbf{x} \\ \sqrt{\lambda} \cdot \mathrm{imag}(\mathbf{B}\mathbf{x}) \end{bmatrix} - \begin{bmatrix} \mathbf{b} \\ \mathbf{0} \end{bmatrix}  \right\|_2^2,
\end{split}\label{eq:psense}
\end{equation}
where $\lambda \in \mathbb{R}$ is a positive regularization parameter and the matrix $\mathbf{B}$ embeds prior information about the image phase such that the regularization encourages $\mathbf{B}\mathbf{x}$ to be real-valued.  Another line of research within MRI instead imposes phase constraints by leveraging linear predictability and the conjugate-symmetry characteristics of the Fourier transform, leading to an inverse problem formulation that can take the general form \cite{haldar2014, haldar2020,haldar2015b,kim2018}
\begin{equation}
\begin{split}
\hat{\mathbf{x}} &= \argmin_{\mathbf{x} \in \mathbb{C}^N} \|\mathbf{A}\mathbf{x} - \mathbf{b}\|_2^2 + \lambda \|\mathbf{C}\mathbf{x} - \mathbf{D} \overline{(\mathbf{E} \mathbf{x})}\|_2^2\\
&= \argmin_{\mathbf{x} \in \mathbb{C}^N} \left\|\begin{bmatrix} \mathbf{A}\mathbf{x} \\ \sqrt{\lambda}\mathbf{C}\mathbf{x} - \sqrt{\lambda}\mathbf{D} \overline{(\mathbf{E} \mathbf{x})} \end{bmatrix} - \begin{bmatrix} \mathbf{b} \\ \mathbf{0} \end{bmatrix}  \right\|_2^2,
\end{split}\label{eq:loraks}
\end{equation}
for appropriate matrices $\mathbf{C}$, $\mathbf{D}$, and $\mathbf{E}$.

Although these are nonlinear least-squares problems because the operators involved are nonlinear, previous work has benefitted from the fact that this kind of inverse problem can be transformed into an equivalent higher-dimensional real-valued linear least-squares problem \cite{bydder2005, willig2005, lew2007, hoge2007, haldar2012,blaimer2016,haldar2014, haldar2020,haldar2015b,kim2018}. Specifically, this can be done by replacing all complex-valued quantities with real-valued quantities, e.g., separating $\mathbf{x} \in \mathbb{C}^N$ into its real and imaginary components,  and treating this as an inverse problem in $\mathbb{R}^{2N}$ rather than the original space $\mathbb{C}^N$.   While this real-valued transformation of the problem is effective and enables the use of standard linear least-squares solution methods, it can also cause computational inefficiencies and can sometimes be difficult to implement when the operators involved have complicated structure.  

In this work, we describe theory that enables provably-convergent linear least-squares iterative algorithms to be applied to this nonlinear least-squares problem setting, without requiring a real-valued transformation of the original complex-valued vectors and operators.  This can enable both improved computation speed and simplified algorithm implementations.

\section{Background}

\subsection{Linear, Antilinear, and Real-Linear Operators}

In this section, we briefly summarize some definitions and properties of linear and antilinear operators, with simplifications corresponding to our finite-dimensional problem context.  Readers interested in a more detailed and more general treatment are referred to Refs.~\cite{rudin1991, huhtanen2011}.

\begin{definition}[Linear Operator] An operator $\mathcal{F}(\cdot): \mathbb{C}^N \rightarrow \mathbb{C}^M$ is said to be linear (or complex-linear) if it satisfies both additivity 
	\begin{equation}
	\mathcal{F}(\mathbf{x} + \mathbf{y}) = \mathcal{F}(\mathbf{x}) + \mathcal{F}(\mathbf{y}) \text{ for } \forall \mathbf{x},\mathbf{y} \in \mathbb{C}^N
	\end{equation}
	and homogeneity 
	\begin{equation}
	\mathcal{F}(\alpha \mathbf{x} ) = \alpha \mathcal{F}(\mathbf{x}) \text{ for } \forall \mathbf{x} \in \mathbb{C}^N, \forall \alpha \in \mathbb{C}.  
	\end{equation}
	
\end{definition}

\begin{property}
	For any linear operator  $\mathcal{F}(\cdot): \mathbb{C}^N \rightarrow \mathbb{C}^M$, there is a unique matrix $\mathbf{F} \in \mathbb{C}^{M \times N}$ such that $\mathcal{F}(\mathbf{x}) = \mathbf{F}\mathbf{x}$ for $\forall \mathbf{x} \in \mathbb{C}^N$.   
\end{property}

\begin{definition}[Antilinear Operator] An operator $\mathcal{G}(\cdot): \mathbb{C}^N \rightarrow \mathbb{C}^M$ is said to be antilinear (or conjugate-linear) if it satisfies both additivity 
	\begin{equation}
	\mathcal{G}(\mathbf{x} + \mathbf{y}) = \mathcal{G}(\mathbf{x}) + \mathcal{G}(\mathbf{y}) \text{ for } \forall \mathbf{x},\mathbf{y} \in \mathbb{C}^N
	\end{equation}
	and conjugate homogeneity
	\begin{equation}
	\mathcal{G}(\alpha\mathbf{x} ) = \overline{\alpha} \mathcal{G}(\mathbf{x})\text{ for } \forall\mathbf{x} \in \mathbb{C}^N, \forall\alpha \in\mathbb{C}.
	\end{equation}
\end{definition}

\begin{property} For any antilinear operator $\mathcal{G}(\cdot): \mathbb{C}^N \rightarrow \mathbb{C}^M$, there is a unique matrix $\mathbf{G} \in \mathbb{C}^{M \times N}$ such that $\mathcal{G}(\mathbf{x}) = \overline{(\mathbf{G}\mathbf{x})}$ for $\forall \mathbf{x} \in \mathbb{C}^N$.  
\end{property}

Note that by taking the matrix $\mathbf{G}$ as the identity matrix, we observe that applying complex conjugation $\overline{\mathbf{x}}$ is an antilinear operation on the vector $\mathbf{x}$.

\begin{definition}[Real-Linear Operator] An operator $\mathcal{A}(\cdot): \mathbb{C}^N \rightarrow \mathbb{C}^M$ is said to be real-linear if it satisfies both additivity 
	\begin{equation}
	\mathcal{A}(\mathbf{x} + \mathbf{y}) = \mathcal{A}(\mathbf{x}) + \mathcal{A}(\mathbf{y}) \text{ for } \forall \mathbf{x},\mathbf{y} \in \mathbb{C}^N
	\end{equation}
	and homogeneity with respect to real-valued scalars
	\begin{equation}
	\mathcal{A}(\alpha \mathbf{x} ) = \alpha \mathcal{A}(\mathbf{x}) \text{ for } \forall \mathbf{x} \in \mathbb{C}^N, \forall \alpha \in \mathbb{R}.  
	\end{equation}
\end{definition}

Real-linearity is a generalization of both linearity and antilinearity, as can be seen from the following property.

\begin{property} 
	Every real-linear operator $\mathcal{A}(\cdot): \mathbb{C}^N \rightarrow \mathbb{C}^M$ can be uniquely decomposed as the sum of a linear operator and an antilinear operator.  In particular, $\mathcal{A}(\mathbf{x}) = \mathcal{F}(\mathbf{x}) + \mathbf{G}(\mathbf{x})$ for $\forall \mathbf{x} \in \mathbb{C}^N$, where $\mathcal{F}(\cdot): \mathbb{C}^N\rightarrow \mathbb{C}^M$ is the linear operator defined by
	\begin{equation}
	\mathcal{F}(\mathbf{x}) \triangleq \frac{1}{2} \mathcal{A}(\mathbf{x}) - \frac{i}{2} \mathcal{A}(i\mathbf{x})
	\end{equation} 
	and $\mathcal{G}(\cdot): \mathbb{C}^N\rightarrow \mathbb{C}^M$ is the antilinear operator defined by
	\begin{equation}
	\mathcal{G}(\mathbf{x}) \triangleq \frac{1}{2} \mathcal{A}(\mathbf{x}) + \frac{i}{2} \mathcal{A}(i\mathbf{x}).
	\end{equation} 
\end{property}

\begin{property} For any real-linear operator $\mathcal{A}(\cdot): \mathbb{C}^N \rightarrow \mathbb{C}^M$, there  are unique matrices $\mathbf{F},\mathbf{G} \in \mathbb{C}^{M \times N}$ such that $\mathcal{A}(\mathbf{x}) = \mathbf{F}\mathbf{x} + \overline{(\mathbf{G}\mathbf{x})}$ for $\forall \mathbf{x} \in \mathbb{C}^N$.  
\end{property}

Notably, both the $\mathrm{real}(\cdot)$ and $\mathrm{imag}(\cdot)$ operators from Eqs.~\eqref{eq:re} and \eqref{eq:im} are observed to have real-linear form.

\begin{property}
	For any two real-linear operators  $\mathcal{A}_1(\cdot): \mathbb{C}^N \rightarrow \mathbb{C}^M$ and $\mathcal{A}_2(\cdot): \mathbb{C}^N \rightarrow \mathbb{C}^M$, their sum $\mathcal{A}_1(\cdot) + \mathcal{A}_2(\cdot)$ is also a real-linear operator.\label{prop:sum}
\end{property}

\begin{property}
	For any two real-linear operators  $\mathcal{A}_1(\cdot): \mathbb{C}^N \rightarrow \mathbb{C}^P$ and $\mathcal{A}_2(\cdot): \mathbb{C}^P \rightarrow \mathbb{C}^M$, their composition $\mathcal{A}_2(\cdot) \circ \mathcal{A}_1(\cdot): \mathbb{C}^N \rightarrow \mathbb{C}^M \triangleq \mathcal{A}_2(\mathcal{A}_1(\cdot)) )$ is also a real-linear operator.\label{prop:comp}
\end{property}

As can be seen, any operator that can be represented as the summation and/or composition of some operators that are linear and some operators that are antilinear can be viewed as a real-linear operator.  As a result, the scenarios of interest in this paper all involve real-linear operators, and the remainder of this paper will assume that $\mathcal{A}(\cdot)$ obeys real-linearity, and has been decomposed in matrix form as $\mathcal{A}(\mathbf{x}) = \mathbf{F}\mathbf{x} + \overline{(\mathbf{G}\mathbf{x})}$.

\subsection{Real-Valued Transformation of Complex-Valued Least Squares}\label{sec:lin}
Assuming $\mathcal{A}(\cdot)$ is real-linear  as described in the previous subsection, Eq.~\eqref{eq:ls} can be rewritten as
\begin{equation}
\hat{\mathbf{x}} = \argmin_{\mathbf{x} \in \mathbb{C}^N} \| \mathbf{F}\mathbf{x} + \overline{(\mathbf{G}\mathbf{x})} - \mathbf{b} \|_2^2,\label{eq:simp}
\end{equation}
which is a nonlinear least squares problem.   However, as stated in the introduction, previous work  \cite{bydder2005, willig2005, lew2007, hoge2007, haldar2012,blaimer2016,haldar2014, haldar2020,haldar2015b,kim2018} has transformed this problem into the form of a conventional linear least-squares problem by treating the variable $\mathbf{x}$ as an element of $\mathbb{R}^{2N}$ instead of $\mathbb{C}^N$.   This was achieved by rewriting $\mathbf{x} \in \mathbb{C}^N$ as $\mathbf{x} = \mathbf{x}_r + i \mathbf{x}_i$, where the real-valued vectors $\mathbf{x}_r, \mathbf{x}_i \in \mathbb{R}^N$ represent the real and imaginary components of $\mathbf{x}$.  This allows us to equivalently rewrite the solution to Eq.~\eqref{eq:simp} as $\hat{\mathbf{x}} = \hat{\mathbf{x}}_r + i \hat{\mathbf{x}}_i$, with
\begin{equation}
\begin{split}
\{\hat{\mathbf{x}}_r,\hat{\mathbf{x}}_i\} &= \argmin_{\mathbf{x}_r,\mathbf{x}_i \in \mathbb{R}^N}\| \mathbf{F}\mathbf{x}_r + i\mathbf{F}\mathbf{x}_i + \overline{\mathbf{G}} \mathbf{x}_r -i \overline{\mathbf{G}} \mathbf{x}_i - \mathbf{b} \|_2^2\\
&= \argmin_{\mathbf{x}_r,\mathbf{x}_i \in \mathbb{R}^N}\left\| \begin{bmatrix} \mathrm{real}(\mathbf{F}\mathbf{x}_r + i\mathbf{F}\mathbf{x}_i + \overline{\mathbf{G}} \mathbf{x}_r -i \overline{\mathbf{G}} \mathbf{x}_i - \mathbf{b}) \\ \mathrm{imag}(\mathbf{F}\mathbf{x}_r + i\mathbf{F}\mathbf{x}_i + \overline{\mathbf{G}} \mathbf{x}_r -i \overline{\mathbf{G}} \mathbf{x}_i - \mathbf{b}) \end{bmatrix} \right\|_2^2\\
&= \argmin_{\tilde{\mathbf{x}} \in \mathbb{R}^{2N}} \left\| \tilde{\mathbf{A}} \tilde{\mathbf{x}} -  \tilde{\mathbf{b}} \right\|_2^2,
\end{split}\label{eq:red}
\end{equation}
where
\begin{equation}
\tilde{\mathbf{x}} \triangleq \begin{bmatrix} \mathbf{x}_r \\ \mathbf{x}_i \end{bmatrix} \in \mathbb{R}^{2N},
\end{equation}
\begin{equation}
\tilde{\mathbf{A}} \triangleq \begin{bmatrix} \mathrm{real}(\mathbf{F}) + \mathrm{real}(\mathbf{G}) & - \mathrm{imag}(\mathbf{F}) - \mathrm{imag}(\mathbf{G}) \\ \mathrm{imag}(\mathbf{F}) - \mathrm{imag}(\mathbf{G}) & \mathrm{real}(\mathbf{F}) - \mathrm{real}(\mathbf{G}) \end{bmatrix}  \in \mathbb{R}^{2M \times 2N},\label{eq:atilde}
\end{equation}
and
\begin{equation}
\tilde{\mathbf{b}} \triangleq \begin{bmatrix} \mathrm{real}(\mathbf{b}) \\ \mathrm{imag}(\mathbf{b}) \end{bmatrix}  \in \mathbb{R}^{2M}.
\end{equation}
The final expression in Eq.~\eqref{eq:red} has the form of a standard real-valued linear least-squares problem, and therefore can be solved using any of the linear least-squares solution methods described in the introduction.  For example,  the Landweber iteration \cite{landweber1951} applied to this problem would proceed as given in Algorithm~\ref{alg:land}, and with infinite numerical precision, $\hat{\mathbf{x}}_k$ is guaranteed to converge to a globally optimal solution as $k \rightarrow \infty$ whenever $0 < \alpha < 2/ \|\tilde{\mathbf{A}}\|_2^2$.

\begin{algorithm}[ht]
	\caption{Landweber Iteration applied to Eq.~\eqref{eq:red}}
	{\bf Inputs:} $\tilde{\mathbf{A}} \in \mathbb{R}^{2M \times 2N}$, $\tilde{\mathbf{b}}\in \mathbb{R}^{2M}$,  $\tilde{\mathbf{x}}_0\in \mathbb{R}^{2N}$ (initial guess for $\tilde{\mathbf{x}}$), and $\alpha \in \mathbb{R}$ \\ \text{\hspace{0.5in}} (step size parameter)\\
	{\bf Initialization:}\\
	\text{\hspace{0.3in}} $k = 0$; \\
	{\bf Iteration:}\\
	\text{\hspace{0.3in}} While stopping conditions are not met:\\
	\text{\hspace{0.6in}} $\tilde{\mathbf{x}}_{k+1}  = \tilde{\mathbf{x}}_{k} + \alpha \tilde{\mathbf{A}}^H (\tilde{\mathbf{b}}- \tilde{\mathbf{A}}\tilde{\mathbf{x}}_k)$; \\  
	\text{\hspace{0.6in}} $k = k+1$;\\
	{\bf Output:} Final value of $\tilde{\mathbf{x}}_{k+1}$
	\label{alg:land}
\end{algorithm}

As another example, the CG algorithm \cite{hestenes1952} applied to this problem would proceed as given in Algorithm~\ref{alg:cg}, and with infinite numerical precision, $\hat{\mathbf{x}}_k$ would be guaranteed to converge to a globally optimal solution after at most $2N$ iterations.  

\begin{algorithm}[ht]
	\caption{Conjugate Gradient Algorithm applied to Eq.~\eqref{eq:red}}
	
	{\bf Inputs:} $\tilde{\mathbf{A}} \in \mathbb{R}^{2M \times 2N}$, $\tilde{\mathbf{b}}\in \mathbb{R}^{2M}$, and $\tilde{\mathbf{x}}_0\in \mathbb{R}^{2N}$ (initial guess for $\tilde{\mathbf{x}}$)\\
	{\bf Initialization:}\\
	\text{\hspace{0.3in}} $\mathbf{r}_0 = \tilde{\mathbf{A}}^H(\tilde{\mathbf{b}} - \tilde{\mathbf{A}} \tilde{\mathbf{x}}_0$);\\
	\text{\hspace{0.3in}} $\mathbf{p}_0 = \mathbf{r}_0$;\\
	\text{\hspace{0.3in}} $k = 0$; \\
	{\bf Iteration:}\\
	\text{\hspace{0.3in}} While stopping conditions are not met:\\
	\text{\hspace{0.6in}} $\mathbf{z}_k = \tilde{\mathbf{A}}^H\tilde{\mathbf{A}} \mathbf{p}_k$; \\ 
	\text{\hspace{0.6in}} $\alpha_k = (\mathbf{r}_k^H\mathbf{r}_k)/(\mathbf{p}_k^H  \mathbf{z}_k)$; \\ 
	\text{\hspace{0.6in}} $\tilde{\mathbf{x}}_{k+1}  = \tilde{\mathbf{x}}_{k} + \alpha_k \mathbf{p}_k$; \\ 
	\text{\hspace{0.6in}} ${\mathbf{r}}_{k+1}  = {\mathbf{r}}_{k} - \alpha_k  \mathbf{z}_k$; \\
	\text{\hspace{0.6in}} $\beta_k = (\mathbf{r}_{k+1}^H\mathbf{r}_{k+1})/(\mathbf{r}_k^H  \mathbf{r}_k)$; \\ 
	\text{\hspace{0.6in}} ${\mathbf{p}}_{k+1}  = {\mathbf{r}}_{k+1} + \beta_k \mathbf{p}_k$; \\ 
	\text{\hspace{0.6in}} $k = k+1$;\\
	{\bf Output:} Final value of $\tilde{\mathbf{x}}_{k+1}$
	\label{alg:cg}
\end{algorithm}
Compared to the analytic linear least-squares solution corresponding to Eq.~\eqref{eq:norm}, these iterative algorithms are generally useful for larger-scale problems where the matrix $\tilde{\mathbf{A}}$ may be too large to store in memory, and where the matrix has structure so that matrix-vector multiplications with $\tilde{\mathbf{A}}$ and $\tilde{\mathbf{A}}^H$ can be computed quickly using specially-coded function calls rather than working with actual matrix representations (e.g., if $\tilde{\mathbf{A}}$ has convolution structure so that matrix-vector multiplication can be implemented using the Fast Fourier Transform, if  $\tilde{\mathbf{A}}$ is sparse, etc.).

Although the problem transformation from Eq.~\eqref{eq:red} has been widely used \cite{bydder2005, willig2005, lew2007, hoge2007, haldar2012,blaimer2016,haldar2014, haldar2020,haldar2015b,kim2018}, it can also be cumbersome to work with if the operator $\mathcal{A}(\cdot)$ has more complicated structure.    For example, the optimization problem in Eq.~\eqref{eq:loraks} involves the composition of linear and antilinear operators, and the $\tilde{\mathbf{A}}$ matrix corresponding to this case has a complicated structure that is laborious to derive.  In particular, with much manipulation, the matrix for this case can be derived to be 
\begin{equation}
\tilde{\mathbf{A}} = \begin{bmatrix} 
\mathrm{real}(\mathbf{A}) & -\mathrm{imag}(\mathbf{A}) \\ 
\mathbf{H}_{11} & \mathbf{H}_{12}   \\ 
\mathrm{imag}(\mathbf{A}) & \mathrm{real}(\mathbf{A})  \\ 
\mathbf{H}_{21} & \mathbf{H}_{22}  \end{bmatrix},\label{eq:loraksbad}
\end{equation}
with
\begin{equation}
\mathbf{H}_{11} = \sqrt{\lambda}\cdot \mathrm{real}(\mathbf{C})- \sqrt{\lambda}\cdot\mathrm{real}(\mathbf{D})\mathrm{real}(\mathbf{E}) - \sqrt{\lambda}\cdot\mathrm{imag}(\mathbf{D})\mathrm{imag}(\mathbf{E}) ,
\end{equation}
\begin{equation}
\mathbf{H}_{12} = -\sqrt{\lambda}\cdot\mathrm{imag}(\mathbf{C}) + \sqrt{\lambda}\cdot\mathrm{real}(\mathbf{D})\mathrm{imag}(\mathbf{E}) - \sqrt{\lambda}\cdot \mathrm{imag}(\mathbf{D})\mathrm{real}(\mathbf{E}) ,
\end{equation}
\begin{equation}
\mathbf{H}_{21} = \sqrt{\lambda}\cdot\mathrm{imag}(\mathbf{C}) - \sqrt{\lambda}\cdot\mathrm{imag}(\mathbf{D})\mathrm{real}(\mathbf{E}) + \sqrt{\lambda}\cdot \mathrm{real}(\mathbf{D})\mathrm{imag}(\mathbf{E}) ,
\end{equation}
and
\begin{equation}
\mathbf{H}_{22} = \sqrt{\lambda}\cdot\mathrm{real}(\mathbf{C}) + \sqrt{\lambda}\cdot\mathrm{imag}(\mathbf{D})\mathrm{imag}(\mathbf{E}) + \sqrt{\lambda}\cdot \mathrm{real}(\mathbf{D})\mathrm{real}(\mathbf{E}).
\end{equation}
Of course, Eq.~\eqref{eq:loraks} relies on a relatively simple mixture of linear and antilinear operators, and problems involving more complicated mixtures would be even more laborious to derive.

Beyond just the effort required to compute the general form of $\tilde{\mathbf{A}}$, it can also be computationally expensive to try to use this type of expression in an iterative algorithm, particularly when the different operators have been implemented as specially-coded function calls.  For example, if we were not given the actual matrix representations of $\mathbf{A}$, $\mathbf{C}$, $\mathbf{D}$, and $\mathbf{E}$ in Eq.~\eqref{eq:loraksbad} and only had function calls that implemented matrix-vector multiplication with these matrices, then a naive implementation of matrix multiplication between $\tilde{\mathbf{A}}$ and a vector would require 4 calls to the function that computes multiplication with $\mathbf{A}$ (e.g., to compute $\mathrm{real}(\mathbf{A})\mathbf{r}$ for an arbitrary real-valued vector $\mathbf{r} \in \mathbb{R}^N$, we could instead compute the complex-valued matrix-vector multiplication function call to obtain $\mathbf{s} = \mathbf{A}\mathbf{r}$, and then use $\mathrm{real}(\mathbf{A})\mathbf{r} = \mathrm{real}(\mathbf{s})$, with  an analogous approach for computing $\mathrm{imag}(\mathbf{A})\mathbf{t}$ for an arbitrary real-valued vector $\mathbf{t} \in \mathbb{R}^N$), 4 calls to the function that computes multiplication with $\mathbf{C}$, 8 calls to the function that computes multiplication with $\mathbf{D}$, and 8 calls to the function that computes multiplication with $\mathbf{E}$. This relatively large number of function calls represents a substantial increase in computational complexity compared to a standard evaluation of the complex-valued forward model, which would only require the use of one function call for each operator.   Of course, this number of computations is based on a naive implementation, and additional careful manipulations could be used to reduce these numbers of function calls by exploiting redundant computations -- however, this would contribute further to the laborious nature of deriving the form of $\tilde{\mathbf{A}}$.

\section{Main Results}\label{sec:new}
Our main results are given by the following lemmas, which enable the use of the real-valued linear least-squares framework from Sec.~\ref{sec:lin} while relying entirely on complex-valued representations and computations.  

\begin{lemma}
	Consider a real-linear operator $\mathcal{A}(\cdot): \mathbb{C}^N \rightarrow \mathbb{C}^M$, with  corresponding $\tilde{\mathbf{A}}$ matrix as defined in Eq.~\eqref{eq:atilde}.   Also consider arbitrary vectors $\mathbf{m} \in \mathbb{C}^N$ and $\mathbf{n} \in \mathbb{C}^M$, which are decomposed into their real and imaginary components according to $\mathbf{m} = \mathbf{m}_r + i\mathbf{m}_i$ and $\mathbf{n} = \mathbf{n}_r + i \mathbf{n}_i$, with $\mathbf{m}_r,\mathbf{m}_i \in \mathbb{R}^N$ and $\mathbf{n}_r,\mathbf{n}_i \in \mathbb{R}^M$.  Then
	\begin{equation}
	\tilde{\mathbf{A}} \begin{bmatrix} \mathbf{m}_r \\ \mathbf{m}_i \end{bmatrix} = \begin{bmatrix} \mathrm{real}(\mathcal{A}(\mathbf{m})) \\ \mathrm{imag}(\mathcal{A}(\mathbf{m})) \end{bmatrix}\label{eq:real}
	\end{equation}
	and
	\begin{equation}
	\tilde{\mathbf{A}}^H \begin{bmatrix} \mathbf{n}_r \\ \mathbf{n}_i \end{bmatrix} = \begin{bmatrix} \mathrm{real}(\mathcal{A}^*(\mathbf{n})) \\ \mathrm{imag}(\mathcal{A}^*(\mathbf{n})) \end{bmatrix},\label{eq:adj}
	\end{equation}
	with $\mathcal{A}^*(\cdot)$ defined below.
\end{lemma}

\begin{definition}[$\mathcal{A}^*(\cdot)$] 
	Consider a real-linear operator $\mathcal{A}(\cdot): \mathbb{C}^N \rightarrow \mathbb{C}^M$, which is represented for $\forall \mathbf{x} \in \mathbb{C}^N$ as $\mathcal{A}(\mathbf{x}) = \mathbf{F}\mathbf{x} + \overline{(\mathbf{G}\mathbf{x})}$ for some matrices $\mathbf{F},\mathbf{G} \in \mathbb{C}^{M \times N}$.  We define $\mathcal{A}^*(\cdot): \mathbb{C}^M \rightarrow \mathbb{C}^N$ as the mapping $ \mathcal{A}^*(\mathbf{n}) \triangleq \mathbf{F}^H \mathbf{n} + \mathbf{G}^H \overline{\mathbf{n}}$ for $\forall \mathbf{n} \in \mathbb{C}^M$.
\end{definition}

Note that $\mathcal{A}^*(\cdot)$ is also a real-linear operator, and can be equivalently written in real-linear form as $\mathcal{A}^*(\mathbf{n}) \triangleq \mathbf{F}^H \mathbf{n} + \overline{(\mathbf{G}^T \mathbf{n})}$ for $\forall \mathbf{n} \in \mathbb{C}^M$, where $^T$ denotes the transpose operation (without conjugation).  Interestingly, it can also be shown that $\mathcal{A}^*(\cdot)$ matches the definition of the adjoint operator of $\mathcal{A}(\cdot)$ from real-linear operator theory \cite{huhtanen2011}.

\begin{lemma}
	Consider a real-linear operator $\mathcal{A}(\cdot): \mathbb{C}^N \rightarrow \mathbb{C}^M$ that can be written as the composition $\mathcal{A}(\cdot) = \mathcal{A}_2(\cdot) \circ \mathcal{A}_1(\cdot)$ of real-linear  operators $\mathcal{A}_1(\cdot): \mathbb{C}^N \rightarrow \mathbb{C}^P$ and $\mathcal{A}_2(\cdot): \mathbb{C}^P \rightarrow \mathbb{C}^M$.     Then $\mathcal{A}^*(\mathbf{n}) = \mathcal{A}_1^*(\mathcal{A}_2^*(\mathbf{n})))$ for $\forall \mathbf{n} \in \mathbb{C}^M$.\label{lemm:sum}
\end{lemma}

\begin{lemma}
	Consider a real-linear operator $\mathcal{A}(\cdot): \mathbb{C}^N \rightarrow \mathbb{C}^M$ that can be written as the summation $\mathcal{A}(\cdot) = \mathcal{A}_1(\cdot) + \mathcal{A}_2(\cdot)$ of real-linear  operators $\mathcal{A}_1(\cdot): \mathbb{C}^N \rightarrow \mathbb{C}^M$ and $\mathcal{A}_2(\cdot): \mathbb{C}^N \rightarrow \mathbb{C}^M$.     Then $\mathcal{A}^*(\mathbf{n}) = \mathcal{A}_1^*(\mathbf{n})+\mathcal{A}_2^*(\mathbf{n})$ for $\forall \mathbf{n} \in \mathbb{C}^M$.\label{lemm:comp}
\end{lemma}

The proofs of these three lemmas are straightforward, and are given in the appendices.  When combined together, these three lemmas completely eliminate the need to derive or work with the real-valued matrix $\tilde{\mathbf{A}}$ in the context of iterative algorithms, because the effects of multiplication with the real-valued matrices $\tilde{\mathbf{A}}$ and $\tilde{\mathbf{A}}^H$ can be obtained equivalently using the complex-valued nonlinear operators $\mathcal{A}(\cdot)$ and $\mathcal{A}^*(\cdot)$.  This can also lead to computational savings, since e.g., computing $\mathrm{real}(\mathcal{A}(\mathbf{m}))$ and $\mathrm{imag}(\mathcal{A}(\mathbf{m}))$ (as needed for computing multiplication of the matrix $\tilde{\mathbf{A}}$ with a vector using Eq.~\eqref{eq:real}) only requires a single call to the function that computes $\mathcal{A}(\mathbf{m})$.  Likewise, computing multiplication of the matrix $\tilde{\mathbf{A}}^H$ with a vector only requires a single call to the function that computes $\mathcal{A}^*(\cdot)$.  And further, if $\mathcal{A}(\cdot)$ is represented as a complicated summation and/or composition of real-linear operators, we can rely on Properties \ref{prop:sum} and \ref{prop:comp} and Lemmas \ref{lemm:sum} and \ref{lemm:comp} to work incrementally with the individual constituent operators, rather than having to work with the monolithic composite operator in its entirety.

As a consequence of these lemmas, it is, e.g., possible to replace the real-valued Landweber iteration from Algorithm~\ref{alg:land} with the  simpler complex-valued iteration given by Algorithm~\ref{alg:landc}.

\begin{algorithm}[ht]
	\caption{Proposed Complex-Valued Landweber Iteration}
	
	{\bf Inputs:} $\mathcal{A}(\cdot): \mathbb{C}^N \rightarrow \mathbb{C}^N$, $\mathbf{b}\in \mathbb{C}^M$, $\mathbf{x}_0 \in \mathbb{C}^N$ (initial guess for ${\mathbf{x}}$), and $\alpha \in \mathbb{R}$ \\ \text{\hspace{0.5in}} (step size parameter)\\
	{\bf Initialization:}\\
	\text{\hspace{0.3in}} $k = 0$; \\
	{\bf Iteration:}\\
	\text{\hspace{0.3in}} While stopping conditions are not met:\\
	\text{\hspace{0.6in}} ${\mathbf{x}}_{k+1}  = {\mathbf{x}}_{k} + \alpha {\mathcal{A}}^* ({\mathbf{b}}- {\mathcal{A}}{\mathbf{x}}_k)$; \\  
	\text{\hspace{0.6in}} $k = k+1$;\\
	{\bf Output:} Final value of ${\mathbf{x}}_{k+1}$
	\label{alg:landc}
\end{algorithm}
With infinite numerical precision, Algorithm~\ref{alg:landc} will produce the exact same sequence of iterates as Algorithm~\ref{alg:land}, and will therefore have the exact same global convergence guarantees stated previously for Landweber iteration.  

We can make similar modifications to the CG algorithm from Algorithm~\ref{alg:cg}, although need the following additional property to be able to correctly handle the inner-products appearing in the CG algorithm.

\begin{property}
	Consider arbitrary vectors $\mathbf{p},\mathbf{q} \in \mathbb{C}^N$, which are decomposed into their real and imaginary components according to $\mathbf{p} = \mathbf{p}_r + i\mathbf{p}_i$ and $\mathbf{q} = \mathbf{q}_r + i \mathbf{q}_i$, with $\mathbf{p}_r,\mathbf{p}_i,\mathbf{q}_r,\mathbf{q}_i \in \mathbb{R}^N$.  Define $\tilde{\mathbf{p}},\tilde{\mathbf{q}}  \in \mathbb{R}^{2N}$  according to 
	\begin{equation}
	\tilde{\mathbf{p}} = \begin{bmatrix} \mathbf{p}_r \\ \mathbf{p}_i  \end{bmatrix} \text{ and } \tilde{\mathbf{q}} = \begin{bmatrix} \mathbf{q}_r \\ \mathbf{q}_i  \end{bmatrix}
	\end{equation}
	Then $\tilde{\mathbf{p}}^H\tilde{\mathbf{q}} = \mathrm{real}(\mathbf{p}^H\mathbf{q})$.
	
\end{property}

Combining this property with the previous lemmas leads to the simple complex-valued iteration for the CG algorithm given by Algorithm~\ref{alg:cgc}.

\begin{algorithm}[ht]
	\caption{Proposed Complex-Valued Conjugate Gradient Algorithm }
	
	{\bf Inputs:} $\mathcal{A}(\cdot): \mathbb{C}^N \rightarrow \mathbb{C}^N$, $\mathbf{b}\in \mathbb{C}^M$, and $\mathbf{x}_0 \in \mathbb{C}^N$ (initial guess for ${\mathbf{x}}$)\\
	{\bf Initialization:}\\
	\text{\hspace{0.3in}} $\mathbf{r}_0 = {\mathcal{A}}^*({\mathbf{b}} - {\mathcal{A}} ({\mathbf{x}}_0)$);\\
	\text{\hspace{0.3in}} $\mathbf{p}_0 = \mathbf{r}_0$;\\
	\text{\hspace{0.3in}} $k = 0$; \\
	{\bf Iteration:}\\
	\text{\hspace{0.3in}} While stopping conditions are not met:\\
	\text{\hspace{0.6in}} $\mathbf{z}_k = {\mathcal{A}}^*({\mathcal{A}}( \mathbf{p}_k))$; \\ 
	\text{\hspace{0.6in}} $\alpha_k = (\mathbf{r}_k^H\mathbf{r}_k)/\mathrm{real}(\mathbf{p}_k^H  \mathbf{z}_k)$; \\ 
	\text{\hspace{0.6in}} ${\mathbf{x}}_{k+1}  = {\mathbf{x}}_{k} + \alpha_k \mathbf{p}_k$; \\ 
	\text{\hspace{0.6in}} ${\mathbf{r}}_{k+1}  = {\mathbf{r}}_{k} - \alpha_k  \mathbf{z}_k$; \\
	\text{\hspace{0.6in}} $\beta_k = (\mathbf{r}_{k+1}^H\mathbf{r}_{k+1})/(\mathbf{r}_k^H  \mathbf{r}_k)$; \\ 
	\text{\hspace{0.6in}} ${\mathbf{p}}_{k+1}  = {\mathbf{r}}_{k+1} + \beta_k \mathbf{p}_k$; \\ 
	\text{\hspace{0.6in}} $k = k+1$;\\
	{\bf Output:} Final value of ${\mathbf{x}}_{k+1}$
	\label{alg:cgc}
\end{algorithm}

While we have only shown complex-valued adaptations of the Landweber and CG algorithms, this same approach is easily applied to other related algorithms like LSQR \cite{paige1982}.

\section{Useful Relations for Common Real-Linear Operators}

Before demonstrating the empirical characteristics of our proposed new approach, we believe that our proposed framework will be easier to use if we enumerated some of the most  common real-linear $\mathcal{A}(\cdot)$ operators and their corresponding $\mathcal{A}^*(\cdot)$ operators. Such a list is provided in Table~\ref{tab}.

\begin{table}[ht]
	\centering
	\begin{tabular}{l|c|c|c|}
		&  \makecell{$\mathcal{A}(\mathbf{x})$ for \\ $\mathbf{x} \in \mathbb{C}^N$} 
		& \makecell{$\mathcal{A}^*(\mathbf{y})$ for \\ $\mathbf{y} \in \mathbb{C}^M$} 
		& \makecell{$\mathcal{A}^*(\mathcal{A}(\mathbf{x}))$ for \\ $\mathbf{x} \in \mathbb{C}^N$} \\ \hline 
		Real-linear  
		& $\mathbf{F}\mathbf{x} + \overline{(\mathbf{G}\mathbf{x})}$ 
		& $\mathbf{F}^H\mathbf{y} + \mathbf{G}^H \overline{\mathbf{y}}$ 
		&  \\ \hline
		Conjugation 
		& $\overline{\mathbf{x}} $ 
		& $\overline{\mathbf{y}}$ 
		& $\mathbf{x}$ \\ \hline
		Real part 
		& $\mathrm{real}(\mathbf{x})$ 
		& $\mathrm{real}(\mathbf{y})$ 
		& $\mathrm{real}(\mathbf{x})$ \\ \hline 
		Imaginary part 
		& $\mathrm{imag}(\mathbf{x})$ 
		& $i \cdot \mathrm{real}(\mathbf{y})$  
		& $i\cdot \mathrm{imag}(\mathbf{x})$ \\ \hline 
		\makecell[l]{System from  \\ Eq.~\eqref{eq:psense} } & $\begin{bmatrix} \mathbf{A}\mathbf{x} \\ \sqrt{\lambda} \cdot \mathrm{imag}(\mathbf{B}\mathbf{x}) \end{bmatrix}$ 
		& \makecell[l]{$\mathbf{A}^H\mathbf{y}_1$ \\ $\hspace{0.05in}+\sqrt{\lambda}i\mathbf{B}^H\mathrm{real}(\mathbf{y}_2)$} 
		& \makecell[l]{$\mathbf{A}^H\mathbf{A}\mathbf{x}$ \\ $\hspace{0.05in}+ \lambda i\mathbf{B}^H\mathrm{imag}(\mathbf{B}\mathbf{x})$} \\ \hline 
		\makecell[l]{System from  \\ Eq.~\eqref{eq:loraks} } & $\begin{bmatrix} \mathbf{A}\mathbf{x} \\ \sqrt{\lambda}\mathbf{C}\mathbf{x} - \sqrt{\lambda}\mathbf{D} \overline{(\mathbf{E} \mathbf{x})} \end{bmatrix}$ 
		& \makecell[l]{$\mathbf{A}^H\mathbf{y}_1$ \\$\hspace{0.3in}+\sqrt{\lambda} \mathcal{B}^*(\mathbf{y}_2)$ } 
		& \makecell[l]{$\mathbf{A}^H\mathbf{A} \mathbf{x}$ \\$\hspace{0.3in}+ \lambda \mathcal{B}^*(\mathcal{B}(\mathbf{x}))$}\\ \hline 
	\end{tabular}
	\caption{Table of common real-linear $\mathcal{A}(\cdot)$ operators and corresponding $\mathcal{A}^*(\cdot)$ operators. We also provide expressions for $\mathcal{A}^*(\mathcal{A}(\cdot))$ in cases where the combined operator takes a simpler form than applying each operator sequentially.  In the last two rows, it is assumed that the matrix $\mathbf{A} \in \mathbb{C}^{M_1 \times N}$, and that the vector $\mathbf{y} \in \mathbb{C}^M$ is divided into two components $\mathbf{y}_1 \in \mathbb{C}^{M_1}$ and $\mathbf{y}_2 \in \mathbb{C}^{M-M_1}$ with $\mathbf{y} = \protect\begin{bmatrix} \mathbf{y}_1^T & \mathbf{y}_2^T \protect\end{bmatrix}^T$.  In the last row, we take $\mathcal{B}(\mathbf{x}) \triangleq \mathbf{C}\mathbf{x} - \mathbf{D} \overline{(\mathbf{E} \mathbf{x})}$, with corresponding $\mathcal{B}^*(\mathbf{y}) = \mathbf{C}^H\mathbf{y} - \mathbf{E}^H \overline{(\mathbf{D}^H\mathbf{y})} $. Note that a special case of equivalent complex-valued operators associated with Eq.~\eqref{eq:psense} (with $\mathbf{B}$ chosen as the identity matrix) was previously presented by Ref.~\cite{bydder2005}, although without the more general real-linear mathematical framework developed in this work. }\label{tab}
\end{table}

\section{Numerical Example}
To demonstrate the potential benefits of our proposed complex-valued approach, we will consider an instance of the problem described by Eq.~\eqref{eq:loraks}.  In this case, the use of complex-valued operations can lead to both a simpler problem formulation and faster numerical computations.  

To address simplicity, we hope that it is obvious by inspection that the process of deriving $\tilde{\mathbf{A}}$ for this case (as given in Eq.~\eqref{eq:loraksbad}, and needed for the conventional real-valued iterative computations) was non-trivial and labor-intensive, while the derivation of $\mathcal{A}(\cdot)$ and $\mathcal{A}^*(\cdot)$ (as given in Table~\ref{tab}, and needed for the proposed new complex-valued iterative computations) was comparatively fast and easy.  

To address the computational benefits of the proposed approach, we will consider a specific realization of Eq.~\eqref{eq:loraks}, in which $\mathbf{x} \in \mathbb{C}^{1000}$, $\mathbf{n} \in \mathbb{C}^{20000}$, $\mathbf{A} \in \mathbb{C}^{20000\times 1000}$, $\mathbf{C} \in \mathbb{C}^{30000 \times 1000}$, $\mathbf{D} \in \mathbb{C}^{30000 \times 2000}$, and $\mathbf{E} \in \mathbb{C}^{2000 \times 1000}$, with the real and imaginary parts of all of these vectors and matrices drawn at random from the i.i.d. Gaussian distribution.  We then took $\mathbf{b} = \mathbf{A}\mathbf{x} + \mathbf{n}$, and set $\lambda = 10^{-3}$.  For this random problem instance, we find the optimal nonlinear least-squares solution in four distinct ways:
\begin{itemize}
	\item {\bf Conventional Real-Valued Approach with Matrices.}  We assume that  $\mathbf{A}$, $\mathbf{C}$, $\mathbf{D}$, and $\mathbf{E}$ are available to us in matrix form, such that it is straightforward to directly precompute the real-valued matrix $\tilde{\mathbf{A}} \in \mathbb{R}^{100000\times 2000}$ from Eq.~\eqref{eq:loraksbad}.  We then use this precomputed matrix directly in iterative linear least-squares solution algorithms like Landweber iteration, CG, and LSQR.  Although the form of this $\tilde{\mathbf{A}}$ matrix was complicated to derive, multiplications with the precomputed $\tilde{\mathbf{A}}$ and $\tilde{\mathbf{A}}^H$ matrices within each iteration should be very computationally efficient, particularly since we have taken 4 separate complex-valued matrices $\mathbf{A}$, $\mathbf{C}$, $\mathbf{D}$, and $\mathbf{E}$ that were originally specified by a sum total of $1.12\times 10^8$ complex-valued entries ($2.24 \times 10^8$ real numbers), and replaced them with a single real-valued matrix specified by only $2\times 10^8$ real numbers. 
	
	\item {\bf Proposed Complex-Valued Approach with Matrices.}  As in the previous case, we assume that $\mathbf{A}$, $\mathbf{C}$, $\mathbf{D}$, and $\mathbf{E}$ are available to us in matrix form, which allows us to directly form the $\mathbf{F}$ and $\mathbf{G}$ matrices corresponding to the complex-valued real-linear formulation of the problem.  Specifically, $\mathbf{F}$ was formed as 
	\begin{equation}
	\mathbf{F} = \begin{bmatrix} \mathbf{A} \\ \sqrt{\lambda} \mathbf{C} \end{bmatrix}
	\end{equation}
	and $\mathbf{G}$ was formed as
	\begin{equation}
	\mathbf{G} = \begin{bmatrix} \mathbf{0} \\ -\sqrt{\lambda} \, \overline{\mathbf{D}}\mathbf{E} \end{bmatrix}.
	\end{equation}
	We then used these precomputed matrices to evaluate $\mathcal{A}(\cdot)$ and $\mathcal{A}^*(\cdot)$ as needed in our proposed complex-valued iterative algorithms.
	
	\item {\bf Conventional Real-Valued Approach with Function Calls.} We assume that we do not have direct access to the $\mathbf{A}$, $\mathbf{C}$, $\mathbf{D}$, and $\mathbf{E}$ matrices, but are only given blackbox functions that calculate matrix-vector multiplications with these matrices and their conjugate transposes.  As such, we implement matrix-vector multiplication with $\tilde{\mathbf{A}}$ (and similarly for $\tilde{\mathbf{A}}^H$) naively in each iteration of the conventional iterative linear least-squares solution algorithms, using multiple calls to each of these functions as described in Section~\ref{sec:lin}.  This approach is not expected to be computationally efficient given the large number of function calls, although is simpler to implement than more advanced approaches that might be developed to exploit redundant computations within Eq.~\eqref{eq:loraksbad}.

	\item {\bf Proposed Complex-Valued Approach with Function Calls.}  As in the previous case, we assume that we do not have direct access to the $\mathbf{A}$, $\mathbf{C}$, $\mathbf{D}$, and $\mathbf{E}$ matrices, but are only given blackbox functions that calculate matrix-vector multiplications with these matrices and their conjugate transposes.  We implement the proposed complex-valued iterative algorithms using the techniques described in Section~\ref{sec:new}, using the expressions for $\mathcal{A}(\cdot)$ and $\mathcal{A}^*(\cdot)$ given in Table~\ref{tab}.

\end{itemize}
For the sake of reproducible research, Matlab code corresponding to this example is included as supplementary material.

\begin{figure}[t]
	\centering
	\includegraphics[width=1\linewidth]{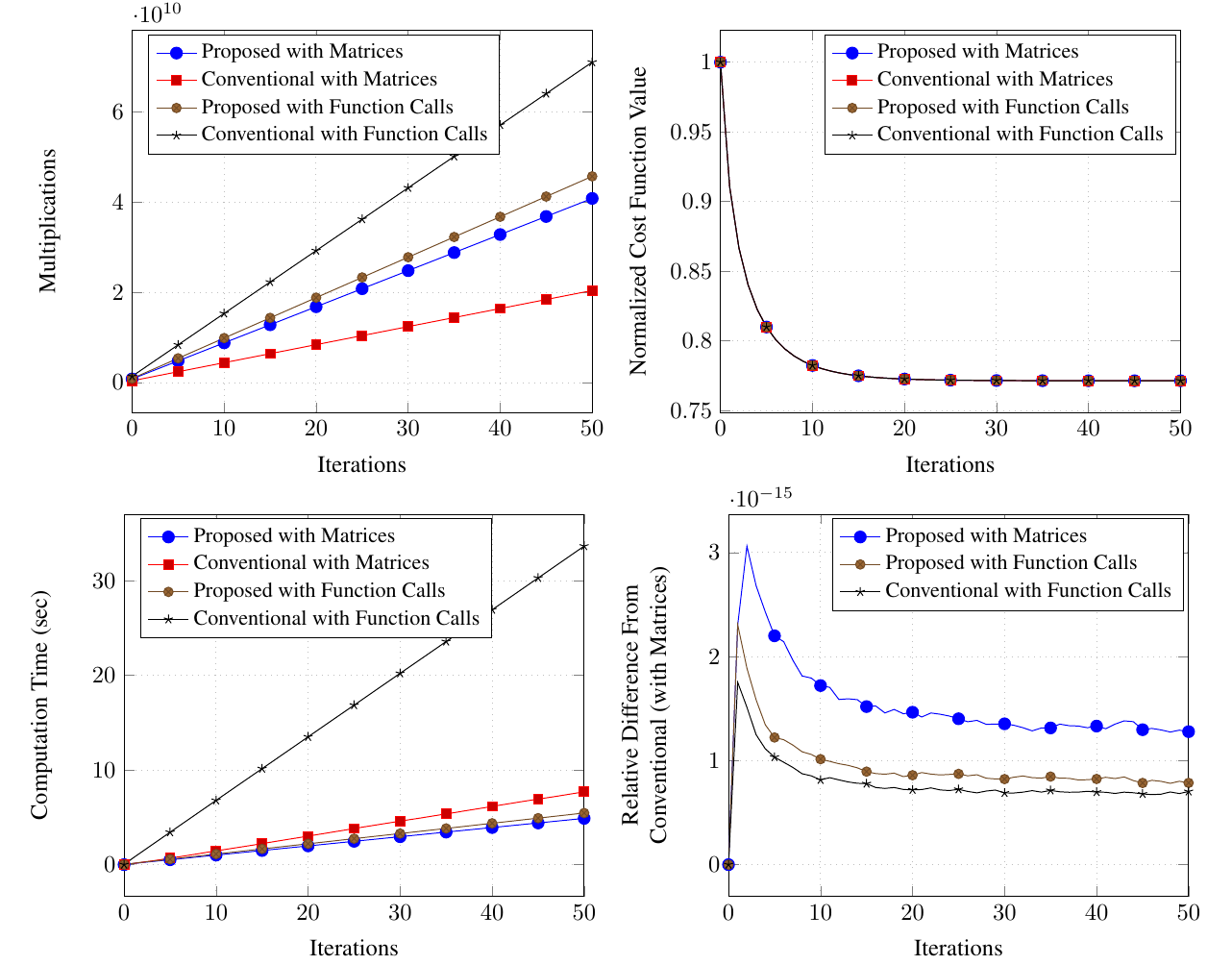} 	
	\caption{ Results for Landweber iteration. The plots show the total number of multiplications, the normalized cost function value (normalized so that the initial value is 1), the computation time in seconds, and the relative difference between the solution from the conventional method with matrices and solutions obtained with other methods.  }
	\label{fig:landweber}
\end{figure}

\begin{figure}[t]
	\centering
	\includegraphics[width=1\linewidth]{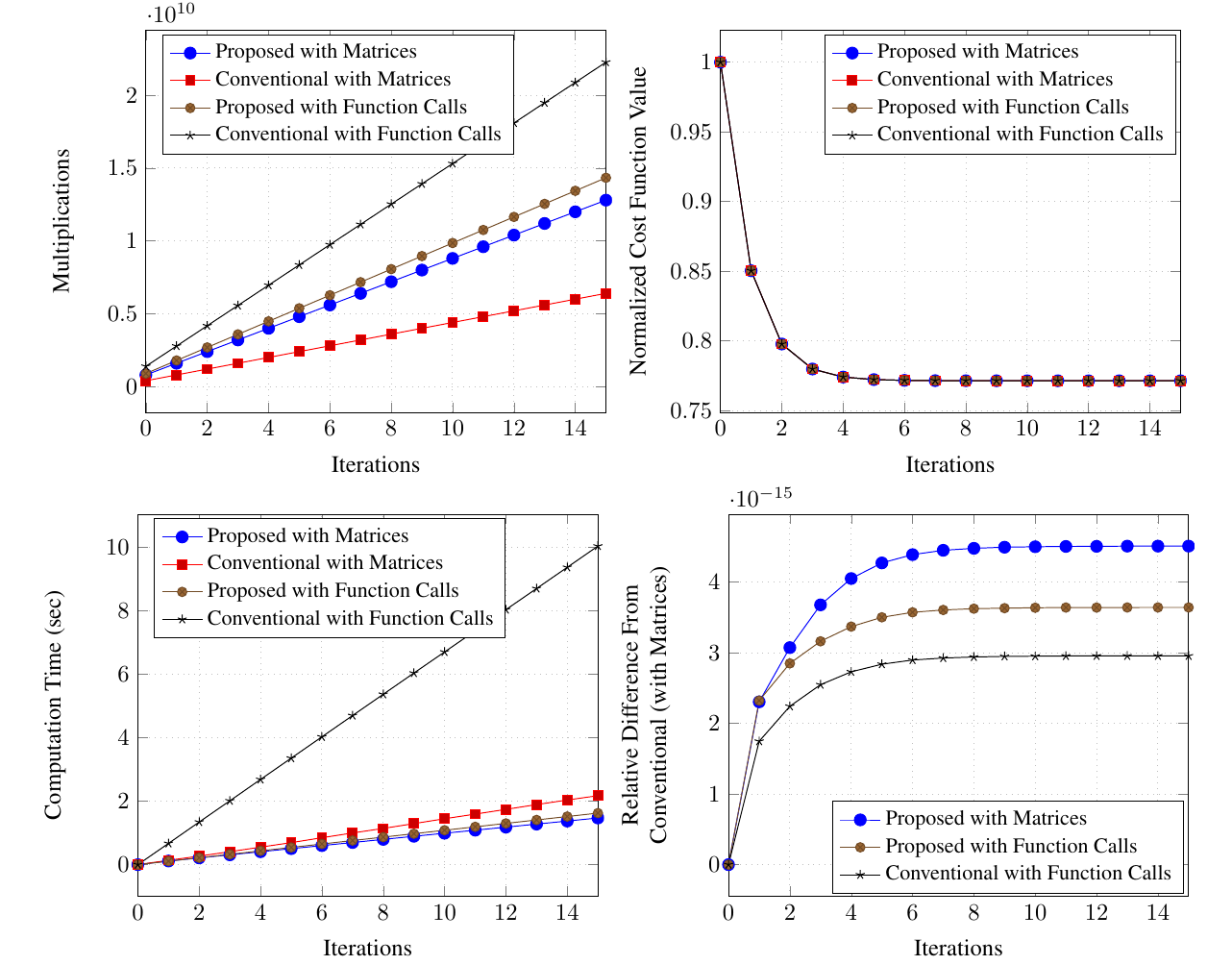} 	
	\caption{Results for the conjugate gradient algorithm. The plots show the total number of multiplications, the normalized cost function value (normalized so that the initial value is 1), the computation time in seconds, and the relative difference between the solution from the conventional method with matrices and solutions obtained with other methods.  }
	\label{fig:conj_grad}
\end{figure}

\begin{figure}[t]
	\centering
	\includegraphics[width=1\linewidth]{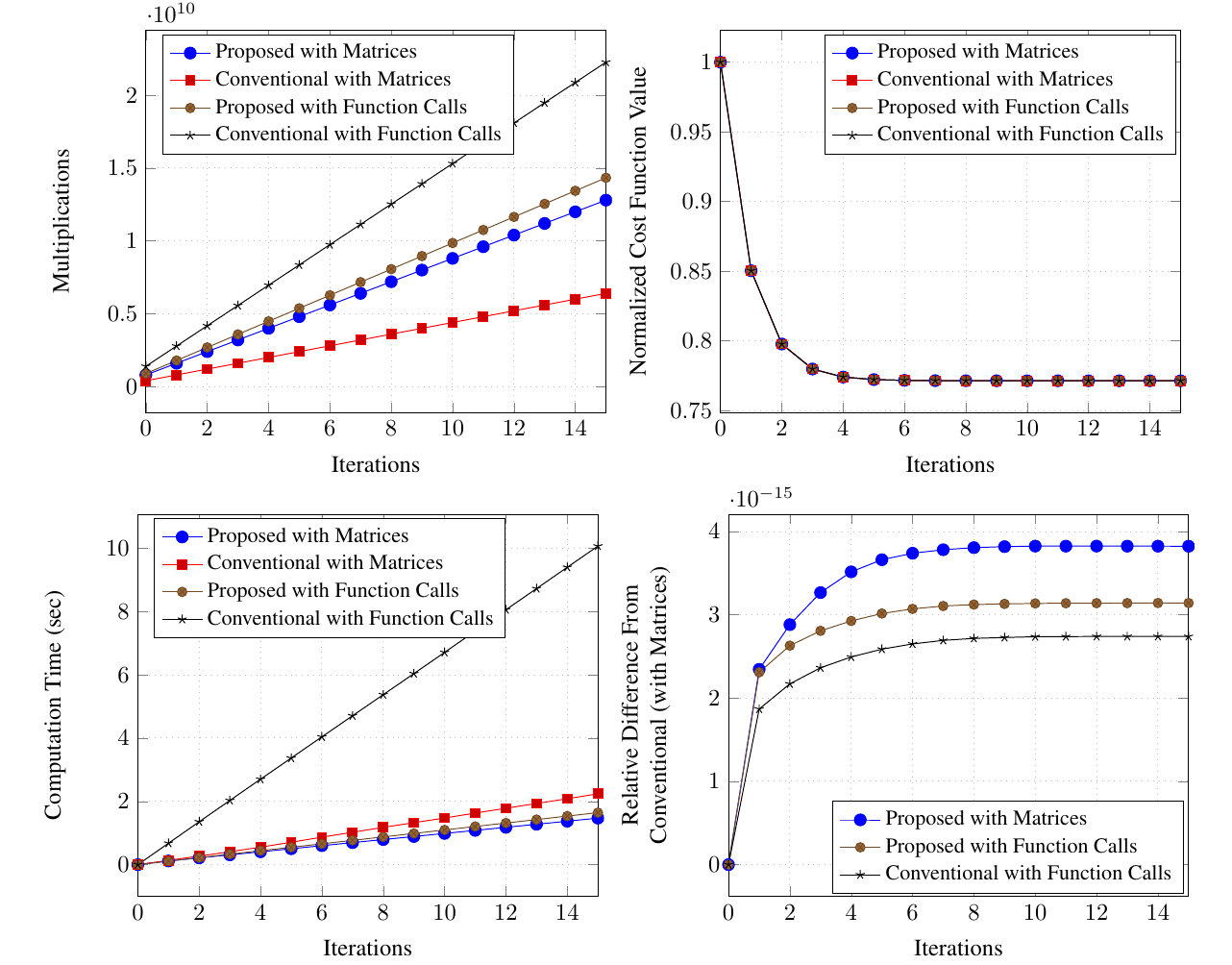} 	
	\caption{Results for the LSQR algorithm. The plots show the total number of multiplications, the normalized cost function value (normalized so that the initial value is 1), the computation time in seconds, and the relative difference between the solution from the conventional method with matrices and solutions obtained with other methods. }
	\label{fig:lsqr}
\end{figure}

For each case, we ran 50 iterations of Landweber iteration and 15 iterations of CG and LSQR in MATLAB 2018b, on a system with an Intel Core i7-8700K 3.70 GHz CPU processor.  For each approach, each algorithm, and at each iteration, we computed (1) the total cumulative number of real-valued scalar multiplications (with 1 complex-valued scalar multiplication equal to 4 real-valued scalar multiplications) used by the algorithm thus far; (2) the cost function value from Eq.~\eqref{eq:loraks} using the current estimate (either $\mathbf{x}_k$ or $\tilde{\mathbf{x}}_k$); (3) the total computation time in seconds; and (4) the relative $\ell_2$-norm difference between the $\mathbf{x}_k$ value estimated from the proposed method with function calls and the other methods, where we define the relative $\ell_2$-norm difference between arbitrary vectors $\mathbf{p}$ and $\mathbf{q}$ as $\|\mathbf{p}-\mathbf{q}\|_2/\|\frac{1}{2}\mathbf{p}+\frac{1}{2}\mathbf{q}\|_2$.  To minimize random fluctuations in computation speed due to background processing, the computation times we report represent the average of 15 different identical trials.

Results for Landweber iteration, the CG algorithm, and LSQR are reported in Figs.~\ref{fig:landweber}-\ref{fig:lsqr}, respectively.  Results confirm that, as should be expected from the theory, all of the different approaches yield virtually identical cost function values and virtually identical solution estimates $\mathbf{x}_k/\tilde{\mathbf{x}}_k$ at each iteration for each of the different algorithms.  There are some very minor differences on the order of $10^{-15}$, which can be attributed to numerical effects resulting from finite-precision arithmetic.  In terms of computational complexity, we observe that the matrix-based approaches are generally associated with fewer multiplications than the implementations that use function calls, which should be expected because the matrix-based approaches were able to precompute simpler consolidated matrix representations that were not available to the function call approaches.  

The proposed approaches required a moderate number of multiplications, somewhat intermediate between the conventional approach with matrices (which had the fewest multiplications) and the conventional approach with function calls (which had the most multiplications).  However, in terms of actual computation time, we observe that the conventional approach with function calls was much slower than any of the other three methods, while the other three methods were all similar to one another.  It is perhaps surprising that the computation times are not directly proportional to the number of multiplications, although this discrepancy is likely related to MATLAB's use of efficient parallelized matrix multiplication libraries.  Importantly, we observe that both variations of the proposed approach are quite fast, and have computation times that are quite similar to the conventional real-valued approach with matrices (which, as we mentioned, was expected to have excellent computational efficiency). There was negligible difference between the computation times assocociated with matrices and function call implementations of the proposed method, which was definitely not the case for the conventional approaches.  And in terms of implementation, the proposed approach with function calls was the easiest to implement, since it didn't require us to derive the forms of any special matrices like $\tilde{\mathbf{A}}$, $\mathbf{F}$, or $\mathbf{G}$, we could just directly work with the individual original matrices $\mathbf{A}$, $\mathbf{C}$, $\mathbf{D}$, and $\mathbf{E}$.

\section{Conclusion}
This work proposed a new approach to solving nonlinear least-squares problems involving real-linear operators.  The new approach allows the use of the original complex-valued operators without transforming them into an unwieldy real-valued form.  Theoretically, the approach enables identical iterative results as the conventional real-valued  transformation, but with much simpler implementation options and potentially much faster computations.  We expect the proposed approach to be valuable for solving general complex-valued nonlinear least-squares problems involving real-linear operators.  Note that the proposed complex-valued approach is  also an integral but previously-undescribed component of the most recent version of an open-source MRI reconstruction software package released by the authors \cite{kim2018}.

\section{Acknowledgments}
This work was supported in part by a USC Annenberg Fellowship, a Kwanjeong Educational Foundation Scholarship, NSF research award CCF-1350563, and NIH research awards R21-EB022951, R01-MH116173, R01-NS074980, R01-NS089212, and R33-CA225400.

\section{References}
\bibliographystyle{ieeetr}
\bibliography{reference}

\appendix

\section{ Proof of Lemma 1}
First, note that Eq.~\eqref{eq:real} is a simple consequence of the derivations shown in Eq.~\eqref{eq:red}.  Thus, the validity of Eq.~\eqref{eq:adj} is the only thing that remains to be proved.  

To see that Eq.~\eqref{eq:adj} is valid, note that
\begin{equation}
\begin{split}
\mathcal{A}^*(\mathbf{n}) &= \mathbf{F}^H\mathbf{n} + \mathbf{G}^H\overline{\mathbf{n}} \\
&= \mathbf{F}^H (\mathbf{n}_r + i \mathbf{n}_i) + \mathbf{G}^H (\mathbf{n}_r - i \mathbf{n}_i)\\
&= \left(\mathrm{real}(\mathbf{F}^H) + i \cdot \mathrm{imag}(\mathbf{F}^H) \right)(\mathbf{n}_r + i \mathbf{n}_i) \\
&\quad + \left(\mathrm{real}(\mathbf{G}^H) + i \cdot \mathrm{imag}(\mathbf{G}^H) \right) (\mathbf{n}_r - i \mathbf{n}_i)\\
&=\left(\mathrm{real}(\mathbf{F}^H)\mathbf{n}_r - \mathrm{imag}(\mathbf{F}^H)\mathbf{n}_i + \mathrm{real}(\mathbf{G}^H)\mathbf{n}_r + \mathrm{imag}(\mathbf{G}^H)\mathbf{n}_i\right) \\
&\quad +i\left(\mathrm{imag}(\mathbf{F}^H)\mathbf{n}_r + \mathrm{real}(\mathbf{F}^H)\mathbf{n}_i + \mathrm{imag}(\mathbf{G}^H)\mathbf{n}_r - \mathrm{real}(\mathbf{G}^H)\mathbf{n}_i\right) \\
&=\left(\mathrm{real}(\mathbf{F})^H\mathbf{n}_r + \mathrm{imag}(\mathbf{F})^H\mathbf{n}_i + \mathrm{real}(\mathbf{G})^H\mathbf{n}_r - \mathrm{imag}(\mathbf{G})^H\mathbf{n}_i\right) \\
&\quad +i\left(-\mathrm{imag}(\mathbf{F})^H\mathbf{n}_r + \mathrm{real}(\mathbf{F})^H\mathbf{n}_i - \mathrm{imag}(\mathbf{G})^H\mathbf{n}_r - \mathrm{real}(\mathbf{G})^H\mathbf{n}_i\right), \\
\end{split}\label{eq:a1}
\end{equation}
where the last line of this expression relies on the fact that $\mathrm{imag}(\mathbf{B}^H) = - \mathrm{imag}(\mathbf{B})^H$ for an arbitrary matrix $\mathbf{B}$.  Equation~\eqref{eq:a1} provides a decomposition of $\mathcal{A}^*(\cdot)$ into its real and imaginary components, and is equivalent to
\begin{equation}
\begin{split}
\begin{bmatrix} \mathrm{real}(\mathcal{A}^*(\mathbf{n})) \\ \mathrm{imag}(\mathcal{A}^*(\mathbf{n})) \end{bmatrix} 
&= \begin{bmatrix} \mathrm{real}(\mathbf{F})^H + \mathrm{real}(\mathbf{G})^H & \mathrm{imag}(\mathbf{F})^H - \mathrm{imag}(\mathbf{G})^H \\ -\mathrm{imag}(\mathbf{F})^H - \mathrm{imag}(\mathbf{G})^H & \mathrm{real}(\mathbf{F})^H - \mathrm{real}(\mathbf{G})^H \end{bmatrix}
\begin{bmatrix} \mathbf{n}_r \\ \mathbf{n}_i \end{bmatrix} \\
&= \begin{bmatrix} \mathrm{real}(\mathbf{F}) + \mathrm{real}(\mathbf{G}) & - \mathrm{imag}(\mathbf{F}) - \mathrm{imag}(\mathbf{G}) \\ \mathrm{imag}(\mathbf{F}) - \mathrm{imag}(\mathbf{G}) & \mathrm{real}(\mathbf{F}) - \mathrm{real}(\mathbf{G}) \end{bmatrix}^H 
\begin{bmatrix} \mathbf{n}_r \\ \mathbf{n}_i \end{bmatrix} \\
&= \tilde{\mathbf{A}}^H \begin{bmatrix} \mathbf{n}_r \\ \mathbf{n}_i \end{bmatrix},
\end{split}
\end{equation}
where the last line comes from the definition of $\tilde{\mathbf{A}}$ in Eq.~\eqref{eq:atilde}.  This proves the validity of  Eq.~\eqref{eq:adj}.
\qed

\section{ Proof of Lemma 2}
Let $\mathcal{A}_1(\cdot): \mathbb{C}^N \rightarrow \mathbb{C}^P$ be a real-linear operator that is represented for $\forall\mathbf{x} \in \mathbb{C}^N$ as $\mathcal{A}_1(\mathbf{x}) = \mathbf{F}_1\mathbf{x} + \overline{(\mathbf{G}_1 \mathbf{x})}$ for some matrices $\mathbf{F}_1,\mathbf{G}_1 \in \mathbb{C}^{P\times N}$, and let $\mathcal{A}_2(\cdot): \mathbb{C}^P \rightarrow \mathbb{C}^M$ be a real-linear operator that is represented for $\forall\mathbf{y} \in \mathbb{C}^P$ as $\mathcal{A}_2(\mathbf{y}) = \mathbf{F}_2\mathbf{y} + \overline{(\mathbf{G}_2 \mathbf{y})}$ for some matrices $\mathbf{F}_2,\mathbf{G}_2 \in \mathbb{C}^{M\times P}$.   Then the composition $\mathcal{A}(\cdot) = \mathcal{A}_2(\cdot) \circ \mathcal{A}_1(\cdot)$ can be expressed for $\forall \mathbf{x} \in \mathbb{C}^N$ as
\begin{equation}
\begin{split}
\mathcal{A}(\mathbf{x}) &= \mathcal{A}_2(\mathcal{A}_1(\mathbf{x})) \\
&= \mathcal{A}_2\left(\mathbf{F}_1\mathbf{x} + \overline{(\mathbf{G}_1\mathbf{x})}\right) \\
&=\mathbf{F}_2\left(\mathbf{F}_1\mathbf{x} + \overline{(\mathbf{G}_1\mathbf{x})}\right) + \overline{\left( \mathbf{G}_2 \left(\mathbf{F}_1\mathbf{x} + \overline{(\mathbf{G}_1\mathbf{x})}\right) \right)} \\
&= (\mathbf{F}_2\mathbf{F}_1 + \overline{\mathbf{G}_2}\mathbf{G}_1)\mathbf{x} + \overline{(\overline{\mathbf{F}_2}\mathbf{G}_1+\mathbf{G}_2\mathbf{F}_1)\mathbf{x}}.
\end{split}
\end{equation}
Thus $\mathcal{A}(\cdot)$ can be written in the real-linear form $\mathcal{A}(\mathbf{x}) = \mathbf{F}\mathbf{x} + \overline{(\mathbf{G}\mathbf{x})}$ for $\forall \mathbf{x} \in \mathbb{C}^N$ with $\mathbf{F} \triangleq \mathbf{F}_2\mathbf{F}_1 + \overline{\mathbf{G}_2}\mathbf{G}_1$ and $\mathbf{G} \triangleq \overline{\mathbf{F}_2}\mathbf{G}_1 + \mathbf{G}_2 \mathbf{F}_1 $.  

By Definition 4, we also have that $\mathcal{A}^*(\mathbf{n}) \triangleq \mathbf{F}^H\mathbf{n} + \mathbf{G}^H \overline{\mathbf{n}}$ for $\forall \mathbf{n} \in \mathbb{C}^M$, $\mathcal{A}_1^*(\mathbf{\mathbf{y}}) \triangleq \mathbf{F}_1^H\mathbf{y} + \mathbf{G}_1^H \overline{\mathbf{y}}$ for $\forall \mathbf{y} \in \mathbb{C}^P$, and $\mathcal{A}_2^*(\mathbf{n}) \triangleq \mathbf{F}_2^H\mathbf{n} + \mathbf{G}_2^H \overline{\mathbf{n}}$ for $\forall \mathbf{n} \in \mathbb{C}^M$.  Thus, we have for $\forall \mathbf{n} \in \mathbb{C}^M$ that
\begin{equation}
\begin{split}
\mathcal{A}_1^*(\mathcal{A}_2^*(\mathbf{n})) &= \mathcal{A}_1^*\left(\mathbf{F}_2^H\mathbf{n} + \mathbf{G}_2^H\overline{\mathbf{n}}\right) \\
&=\mathbf{F}_1^H \left(\mathbf{F}_2^H\mathbf{n} + \mathbf{G}_2^H\overline{\mathbf{n}}\right)  + \mathbf{G}_1^H \overline{\left(\mathbf{F}_2^H\mathbf{n} + \mathbf{G}_2^H\overline{\mathbf{n}}\right)} \\
&= (\mathbf{F}_1^H\mathbf{F}_2^H + \mathbf{G}_1^H \overline{\mathbf{G}_2^H})\mathbf{n} + 
\left(\mathbf{F}_1^H\mathbf{G}_2^H + \mathbf{G}_1^H\overline{\mathbf{F}_2^H}\right)\overline{\mathbf{n}} \\
&= (\mathbf{F}_2\mathbf{F}_1 + \overline{\mathbf{G}_2}\mathbf{G}_1)^H \mathbf{n} + ( \overline{\mathbf{F}_2}\mathbf{G}_1 + \mathbf{G}_2 \mathbf{F}_1 )^H\overline{\mathbf{n}} \\ 
&= \mathbf{F}^H\mathbf{n} + \mathbf{G}^H\overline{\mathbf{n}} \\
&= \mathcal{A}^*(\mathbf{n}),
\end{split}
\end{equation}  
which shows that $\mathcal{A}^*(\mathbf{n}) = \mathcal{A}_1^*(\mathcal{A}_2^*(\mathbf{n}))$ for $\forall \mathbf{n} \in \mathbb{C}^M$ as desired.
\qed

\section{ Proof of Lemma 3}
Let $\mathcal{A}_1(\cdot): \mathbb{C}^N \rightarrow \mathbb{C}^M$ be a real-linear operator that is represented for $\forall\mathbf{x} \in \mathbb{C}^N$ as $\mathcal{A}_1(\mathbf{x}) = \mathbf{F}_1\mathbf{x} + \overline{(\mathbf{G}_1 \mathbf{x})}$ for some matrices $\mathbf{F}_1,\mathbf{G}_1 \in \mathbb{C}^{M\times N}$, and let $\mathcal{A}_2(\cdot): \mathbb{C}^N \rightarrow \mathbb{C}^M$ be a real-linear operator that is represented for $\forall\mathbf{x} \in \mathbb{C}^N$ as $\mathcal{A}_2(\mathbf{x}) = \mathbf{F}_2\mathbf{x} + \overline{(\mathbf{G}_2 \mathbf{x})}$ for some matrices $\mathbf{F}_2,\mathbf{G}_2 \in \mathbb{C}^{M\times N}$.   Then the summation $\mathcal{A}(\cdot) = \mathcal{A}_1(\cdot) + \mathcal{A}_2(\cdot)$ can be expressed for $\forall \mathbf{x} \in \mathbb{C}^N$ as
\begin{equation}
\begin{split}
\mathcal{A}(\mathbf{x})&=\mathcal{A}_1(\mathbf{x})+\mathcal{A}_2(\mathbf{x}) \\
&= \mathbf{F}_1 \mathbf{x} + \overline{(\mathbf{G}_1 \mathbf{x})} + \mathbf{F}_2 \mathbf{x} + \overline{(\mathbf{G}_2 \mathbf{x})} \\
&=\left(\mathbf{F}_1 + \mathbf{F}_2 \right)\mathbf{x} + 
\overline{\left(\mathbf{G}_1 + \mathbf{G}_2\right)\mathbf{x}}.
\end{split}
\end{equation}
Thus $\mathcal{A}(\cdot)$ can be written in the real-linear form $\mathcal{A}(\mathbf{x}) = \mathbf{F}\mathbf{x} + \overline{(\mathbf{G}\mathbf{x})}$ for $\forall \mathbf{x} \in \mathbb{C}^N$ with $\mathbf{F} \triangleq \mathbf{F}_1 +\mathbf{F}_2$ and $\mathbf{G} \triangleq \mathbf{G}_1 + \mathbf{G}_2$.  

By Definition 4, we also have that $\mathcal{A}^*(\mathbf{n}) \triangleq \mathbf{F}^H\mathbf{n} + \mathbf{G}^H \overline{\mathbf{n}}$ for $\forall \mathbf{n} \in \mathbb{C}^M$, $\mathcal{A}_1^*(\mathbf{\mathbf{y}}) \triangleq \mathbf{F}_1^H\mathbf{y} + \mathbf{G}_1^H \overline{\mathbf{y}}$ for $\forall \mathbf{y} \in \mathbb{C}^P$, and $\mathcal{A}_2^*(\mathbf{n}) \triangleq \mathbf{F}_2^H\mathbf{n} + \mathbf{G}_2^H \overline{\mathbf{n}}$ for $\forall \mathbf{n} \in \mathbb{C}^M$.  Thus, we have for $\forall \mathbf{n} \in \mathbb{C}^M$ that
\begin{equation}
\begin{split}
\mathcal{A}_1^*(\mathbf{n}) + \mathcal{A}_2^*(\mathbf{n}) &= \mathbf{F}_1^H \mathbf{n} + \mathbf{G}_1^H\overline{\mathbf{n}} + \mathbf{F}_2^H \mathbf{n} + \mathbf{G}_2^H\overline{\mathbf{n}} \\
& = (\mathbf{F}_1^H + \mathbf{F}_2^H)\mathbf{n} + (\mathbf{G}_1^H+\mathbf{G}_2^H)\overline{\mathbf{n}}\\
& = (\mathbf{F}_1 + \mathbf{F}_2)^H\mathbf{n} + (\mathbf{G}_1+\mathbf{G}_2)^H\overline{\mathbf{n}}\\
& = \mathbf{F}^H \mathbf{n} + \mathbf{G}^H\overline{\mathbf{n}}\\
&= \mathcal{A}^*(\mathbf{n}),
\end{split}
\end{equation}
which shows that $\mathcal{A}^*(\mathbf{n}) = \mathcal{A}_1^*(\mathbf{n}) + \mathcal{A}_2^*(\mathbf{n})$ for $\forall \mathbf{n} \in \mathbb{C}^M$ as desired.
\qed
\end{document}